\documentclass[12pt, twoside, leqno]{article}

\usepackage{amsmath,amsthm}
\usepackage{amssymb}

\usepackage{enumerate}

\usepackage{graphicx}

\newcommand{\pe}{\psi}
\def\d{\mbox{d}} 
\def\ds{\displaystyle} 
\def\e{{\epsilon}} 
\def\eb{\bar{\eta}}  
\def\enorm#1{\|#1\|_2} 
\def\Fp{F^\prime}  
\def\fishpack{{FISHPACK}} 
\def\fortran{{FORTRAN}} 
\def\gmres{{GMRES}} 
\def\gmresm{{\rm GMRES($m$)}} 
\def\Kc{{\cal K}} 
\def\norm#1{\|#1\|} 
\def\wb{{\bar w}} 
\def\zb{{\bar z}} 
\def\half{\hbox{$1\over2$}}
\def\qart{\hbox{$3\over2$}}
\def\E{{\rm E}}
\def\Rex{\mbox{I \hspace{-2.00ex} R}}
\def\KL{Kullback--Leibler }
\def\E{{\rm E}}
\def\Dir{\mbox{Dir}}
\def\N{\mbox{N}}
\def\ga{\mbox{Ga}}
\def\ex{\mbox{Ex}}
\def\DP{\mbox{DP}}
\def\po{\mbox{Po}}
\def\Var{{\rm Var}}
\def\P{\mbox{P}}
\def\d{{\rm d}}
\def\iid{\rm independent and identically distributed }

\def\pr{{\rm pr}}
\def\Un{\mbox{Un}}
\def\data{x^n}
\def\midd{\,|\,}
\def\cF{{\cal F}}
\def\cX{{\cal X}}
\def\cI{{\cal I}}
\def\cD{{\cal D}}
\def\cL{{\cal L}}
\def\cC{{\cal C}}
\def\cH{{\cal H}}
\def\cS{{\cal S}}
\def\cB{{\cal B}}
\def\cM{{\cal M}}
\def\arr{\rightarrow}
\def\be{\mbox{be}}
\def\hatt{\widehat}
\def\sumin{\sum_{i=1}^n}
\def\prodin{\prod_{i=1}^n}
\def\eps{\varepsilon}

\def\bfE{\mbox{\boldmath$E$}}
\def\bfG{\mbox{\boldmath$G$}}

\newtheorem{thm}{Theorem}[section]

\newtheorem{lem}[thm]{Lemma}

\theoremstyle{definition}

\numberwithin{equation}{section}

\begin{document}


\baselineskip=17pt

\title{On an Asymptotic Distribution for the MLE}


\author{Stephen G. Walker\\ \\
Department of Mathematics \\
University of Texas at Austin, USA \\
e-mail: s.g.walker@math.utexas.edu}
       
\date{}       


\maketitle

\renewcommand{\thefootnote}{}



\renewcommand{\thefootnote}{\arabic{footnote}}
\setcounter{footnote}{0}

\begin{abstract}
The paper presents a novel asymptotic distribution for a mle when the log--likelihood is strictly concave in the parameter for all data points; for example, the exponential family. The new asymptotic distribution can be seen as a refinement of the usual normal asymptotic distribution and is comparable to an Edgeworth expansion. However, it is obtained with weaker conditions than even those for asymptotic normality. The same technique is then used to find the exact distribution of the weighted likelihood bootstrap sampler.
\end{abstract}


       















\section{Introduction}  One important aspect of statistical inference is quantifying the uncertainty in statistics; for example, the sampling distribution of the maximum likelihood estimator arising from a model and data.  If approximations are required, it is the asymptotic normal distribution which is often, if not always, used. In this paper we show that if the log--likelihood is strictly concave in the parameter for all data sets, then an improved asymptotic distribution is available. The density estimate has similar properties to a second order Edgeworth expansion, which uses up to three derivtives of the log--likelihood (see \cite{pfanzagl} and \cite{kolassa}); whereas we obtain this using only one derivative. It is the concavity of the log--likelihood which facilitates this. It is also clearly to be seen how to get the asymptotic normal distribution from this new asymptotic distribution.

Consider the family of density functions $f(x;\theta)$, with respect to some dominating measure, which will either be the counting measure or the Lebesgue measure. Here $x\in\mathbb{X}$ and $\theta\in\Theta\subset\mathbb{R}$. We write 
$$l(x;\theta)=-\log f(x;\theta),$$ 
negative the score function,
and assume that $l'(x;\theta)=\partial l(x;\theta)/\partial\theta$ exists for all $\theta$ and $x$, and that $l(x;\theta)$ is strictly convex in $\theta$ for all $x$; i.e. for all $\theta\ne\theta'$ is is that 
$$l(x;\theta)>l(x;\theta')+(\theta-\theta')l'(x;\theta').$$
An example of such is the exponential family, see for example, \cite{Kupperman}; so for some functions $c(x)$ and $t(x)$,
$$f(x;\theta)=c(x)\,\exp\{t(x)\theta-b(\theta)\}\quad\mbox{so}\quad l(x;\theta)=-\log c(x)-t(x)\theta+b(\theta)$$
where $b$ is the normalizing constant and known to be a convex function.

Now let $\theta^*$ be a true parameter value which generates independent and identically distributed data $(X_1,\ldots,X_n)$ from $f(x;\theta^*)$. The maximum likelihood estimator is given by the $\widehat{\theta}$ solving
$$\sum_{i=1}^n l'(X_i;\theta)=0.$$
The paper is about an asymptotic distribution for $\widehat{\theta}$ which will be presented in section 2. Before this we highlight the conditions for asymptotic normality, see for example \cite{Lehmann}, which can be found as follows:

\begin{description}

\item (a) The parameter space $\Theta$ is an open interval.

\item (b) The set $A=\{x:f(x;\theta)>0\}$ does not depend on $\theta$.

\item (c) For all $x\in A$ the density $f(x;\theta)$ is thrice diffentiable with respect to $\theta$ and the third derivative is continuous in $\theta$.

\item (d) The Fisher information $I(\theta)$ satisfies
$$0<I(\theta)=\int \big(l'(x;\theta)\big)^2\,f(x;\theta)\,d x<\infty.$$

\item (e) $\E_\theta\,[ l'(x;\theta)]=0$ and $\E_{\theta}\,[l''(x;\theta)]=I(\theta)$. 

\item (f) For $\theta^*\in\Theta$ there exists $c>0$ and function $M(x)$ such that
$|l'''(x,\theta)|\leq M(x)$ for all $x\in A$ and $|\theta-\theta^*|<c$ and
$\E_{\theta^*}\,[M(x)]<\infty$.

\end{description}

\noindent
As is well known, under these conditions,
\begin{equation}\label{mle}
\sqrt{n}\,\,(\widehat{\theta}-\theta^*)\,\sqrt{I(\theta^*)}\to_d \N\big(0,1\big),
\end{equation}
where $\to_d$ denotes convergence in distribution. This we also write as
\begin{equation}\label{norm}
\widehat{F}_{\widehat{\theta}}^{(N)}(z)=\Phi\left((z-\theta^*)\,\sqrt{n\,I(\theta^*})\right),
\end{equation}
as the estimator of the distribution function for the mle. Here the superscript $N$ refers to the normal approximation.

In section 2  we provide the new asymptotic distributions for the mle under the concave condition. Section 3 then presents three illustrations and section 4 uses the same technique to find the exact distribution for the weighted likelihood bootstrap sampler. Finally, section 5 concludes with some ideas for future work and considers the multivariate case.

\section{An asymptotic distribution for $\widehat{\theta}$}

Let us define
$$D(\theta,\theta^*)=\int l'(x;\theta)\,f(x;\theta^*)\,d x\quad
\mbox{and}\quad V(\theta,\theta^*)=\int \big(l'(x;\theta)\big)^2\,f(x;\theta^*)\,d x.$$
Under Assumption (e) we have that $D(\theta^*,\theta^*)=0$ and $V(\theta^*,\theta^*)=I(\theta^*)$.
Defining 
$$T_n(z)=\frac{1}{n}\,\sum_{i=1}^n l'(X_i;z),$$
and note that $T_n(z)=-S_n(z)/n$, where $S_n(z)$ is the usual score function, the
asymptotic normality of $T_n(z)$, for each $z\in\Theta$, implies
$$A_n(z)=\sqrt{n}\frac{T_n(z)-D(z,\theta^*)}{\sqrt{V(z,\theta^*)-D^2(z,\theta^*)}}\to_d \N(0,1).$$
Our estimator of the distribution of the mle is based on this asymptotic result.
Hence,
$$\widehat{F}_{T_n(z)}(t)=\Phi\left(\sqrt{n}\,\,\frac{t-D(z,\theta^*)}{\sqrt{V(z,\theta^*)-D^2(z,\theta^*)}}\right).$$
See \cite{Lehmann} for further details on asymptotic normality of sample means. 

\noindent
The main result of the paper is the following theorem:

\begin{thm}
Under Assumptions (a), (b),  (d) and (e) combined with $l(x;\theta)$ being strictly convex in $\theta$ for all $x$, it is that
\begin{equation}\label{alt}
\widehat{F}_{\widehat{\theta}}(z)= \Phi\left(\sqrt{n}\,\,\frac{D(z,\theta^*)}{\sqrt{V(z,\theta^*)-D^2(z,\theta^*)}}\right)\end{equation}
is an estimator of the distribution of the mle.
\end{thm}

\begin{proof}
Since $l(x;\theta)$ is strictly convex in $\theta$ for all data sets, we have the observation that
$$\P\left(\widehat{\theta}\leq z\right)=\P\left(T_n(z)\geq 0\right).$$ 
Hence, 
$$\P\left(\widehat{\theta}\leq z\right)=\P\left(A_n(z)\geq -\sqrt{n}\,\,\frac{D(z;\theta^*)}{\sqrt{V(z,\theta^*)-D^2(z,\theta^*)}} \right)$$
and since $A_n(z)\to_d \N(0,1)$,
$$\frac{\P\left(\widehat{\theta}\leq z\right)}
{1-\Phi\left(-\sqrt{n}\,\,\frac{D(z;\theta^*)}{\sqrt{V(z,\theta^*)-D^2(z,\theta^*)}}\right)}\to 1$$
for all $z$.
Consequently, we have and can take the estimator
$$\widehat{\P\left(\widehat{\theta}\leq z\right)}= 1-\Phi\left(-\sqrt{n}\,\,\frac{D(z;\theta^*)}{\sqrt{V(z,\theta^*)-D^2(z,\theta^*)}}\right),$$
which is given by (\ref{alt}), completing the proof.
\end{proof}

\noindent
We can see clearly how to get (\ref{mle}) from (\ref{alt}); requiring the approximations 
$$D(z,\theta^*)\approx D(\theta^*,\theta^*)+(z-\theta^*)\,\frac{\partial D}{\partial z}(\theta^*,\theta^*)=(z-\theta^*)\,I(\theta^*),$$
and $V(z,\theta^*)-D^2(z,\theta^*)\approx V(\theta^*,\theta^*)$, noting 
$\partial D(\theta^*,\theta^*)/\partial z=V(\theta^*,\theta^*)=I(\theta^*)$. This involves some rather loose approximations and suggests the normal approximation should not necessarily work well with $z$ away from $\theta^*$. Indeed we see this phenomenon in an illustration which follows.

Before this we see how (\ref{alt}) is comparable to an Edgeworth expansion. The following standard expansion is to be found in Chapter 16 in \cite{dasgupta};
$$\P\left(\sqrt{n}\,\,(\widehat{\theta}-\theta^*)\,\sqrt{I(\theta^*)}\leq x \right)=\Phi(x)+\phi(x)\,(a+bx^2)/\sqrt{n},$$
where $a$ and $b$ use up to the third derivatives of $l(x;\theta)$, and are based on expectations with respect to $f(x;\theta^*)$.

\begin{lem}
From (\ref{alt}) we obtain
$$\P\left(\sqrt{n}\,\,(\widehat{\theta}-\theta^*)\,\sqrt{I(\theta^*)}\leq  x\right)=\Phi(x)+\half \,c\,\phi(x)\,x^2/\sqrt{n},$$
where 
$$c=\frac{\partial^2 D}{\partial z^2}(\theta^*,\theta^*)\,V(\theta^*,\theta^*)^{-3/2}-\frac{\partial V}{\partial z}(\theta^*,\theta^*)\,V(\theta^*,\theta^*)^{-1}.$$ 
\end{lem}

\begin{proof}
The proof to this uses
$$D\left(\theta^*+\frac{d}{\sqrt{n}},\theta^*\right)=\frac{d}{\sqrt{n}}\,\frac{\partial D}{\partial z}(\theta^*,\theta^*)
+\half\,\frac{d^2}{n}\,\frac{\partial^2 D}{\partial z^2}(\theta^*,\theta^*)+O(n^{-3/2})$$
and
$$V\left(\theta^*+\frac{d}{\sqrt{n}},\theta^*\right)=V(\theta^*,\theta^*)+\frac{d}{\sqrt{n}}\,\frac{\partial V}{\partial z}(\theta^*,\theta^*)+O(n^{-1}).$$
\end{proof}

\noindent
So note we recover an Edgeworth expansion type estimate for the distribtuion of the mle; i.e. (\ref{alt}), using only $D$ and $V$ which themselves only depend on the first derivative of $l(x;\theta)$.

\section{Illustrations} 

\subsection{Exponential family} Consider the exponential family with functions $t(x)$ and $b(\theta)$ so
$$D(\theta,\theta^*)=b'(\theta)-b'(\theta^*)\quad\mbox{and}\quad V(\theta,\theta^*)=b''(\theta^*)+(b'(\theta)-b'(\theta^*))^2.$$ 
Therefore,
$$\widehat{F}_{\widehat{\theta}}(z)= \Phi\left(\sqrt{n}\,\,\frac{b'(z)-b'(\theta^*)}{\sqrt{b''(\theta^*)}}\right).$$
On the other hand, the asymptotic normal distribution is given by
$$\widehat{F}^{(N)}_{\widehat{\theta}}(z)\approx \Phi\left(\sqrt{n}\,\,(z-\theta^*)\,\sqrt{b''(\theta^*)}\right).$$
In particular, suppose 
$f(x;\theta)=\theta\,e^{-x\theta}$, with $x>0$ and $\theta>0$. Then $b(\theta)=-\log\theta$ so $b'(\theta)=-1/\theta$ and $b''(\theta)=1/\theta^2$. 

\begin{center}
\begin{figure}[!htbp]
\begin{center}
\includegraphics[scale=0.5]{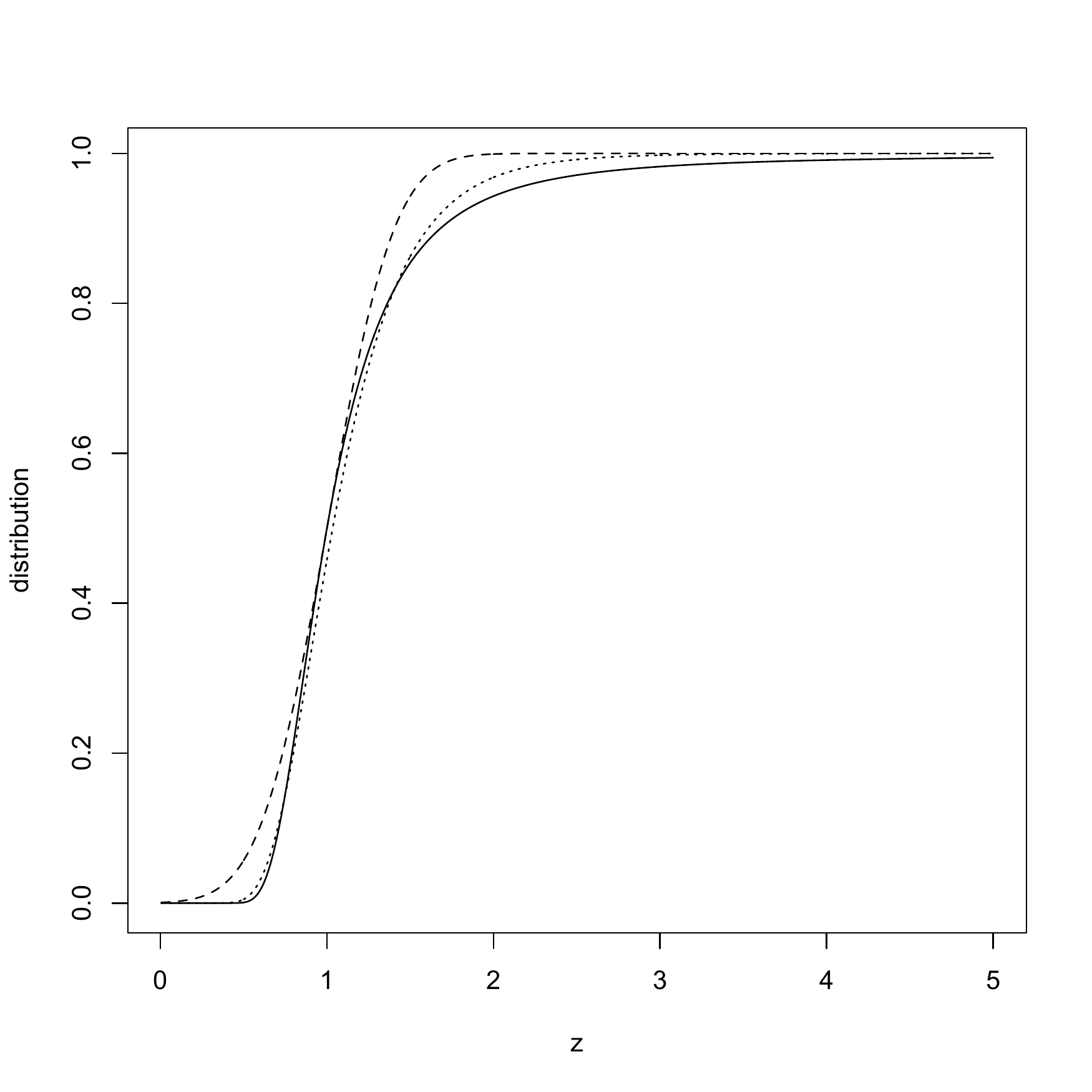}
\caption{(i) Dotted line: $F^*_{\widehat{\theta}}(z)$; (ii) solid line: $\widehat{F}_{\widehat{\theta}}(z)$; (iii) dashed line $\widehat{F}^{(N)}_{\widehat{\theta}}(z)$.}
 \label{fig1}
\end{center}
\end{figure}
\end{center}

Thus, we wish to compare
\begin{equation}\label{fnew}
\widehat{F}_{\widehat{\theta}}(z)=\Phi\left(\sqrt{n}\,\,\theta^{*}\,(1/\theta^*-1/z)\right),
\end{equation}
the new asymptotic distribution, which is not a normal distribution, with
\begin{equation}\label{old}
\widehat{F}^{(N)}_{\widehat{\theta}}(z)=\Phi\left(\sqrt{n}\,\,(z-\theta^*)/\theta^*\right)
\end{equation}
which is the usual asymptotic normal distribution. The true distribution for $\widehat{\theta}$ is
\begin{equation}\label{true}
F^*_{\widehat{\theta}}(z)=1-\Gamma_n(n\theta^*/z),
\end{equation}
where $\Gamma_n$ is the distribution function of a gamma random variable with shape parameter $n$ and scale parameter 1.

With $n=10$ and $\theta^*=1$, Fig.~1 presents three curves; the bold solid line is $\widehat{F}_{\widehat{\theta}}(z)$, the dotted line is the true distribution of $\widehat{\theta}$, whereas the dashed line is 
$\widehat{F}^{(N)}_{\widehat{\theta}}(z)$. So we see that the latter distribution is not accurate whereas the former, even with a sample of size 10 is good.

\subsection{Fisk distribution} The Fisk density function is given by
$$f(x;\theta)=\frac{\theta\,x^{\theta-1}}{(1+x^\theta)^2},\quad x>0$$
and $\theta>0$. Hence
$$l(x;\theta)=2\log(1+x^\theta)-(\theta-1)\,\log x-\log\theta,$$
which is easily shown to be convex in $\theta$ for all $x>0$. Further,
$$l'(x;\theta)=2\frac{\log x\, x^\theta}{1+x^\theta}-\log x-1/\theta.$$
This then gives us access to $D(z,\theta^*)$ and $V(z,\theta^*)$.

\begin{center}
\begin{figure}[!htbp]
\begin{center}
\includegraphics[scale=0.5]{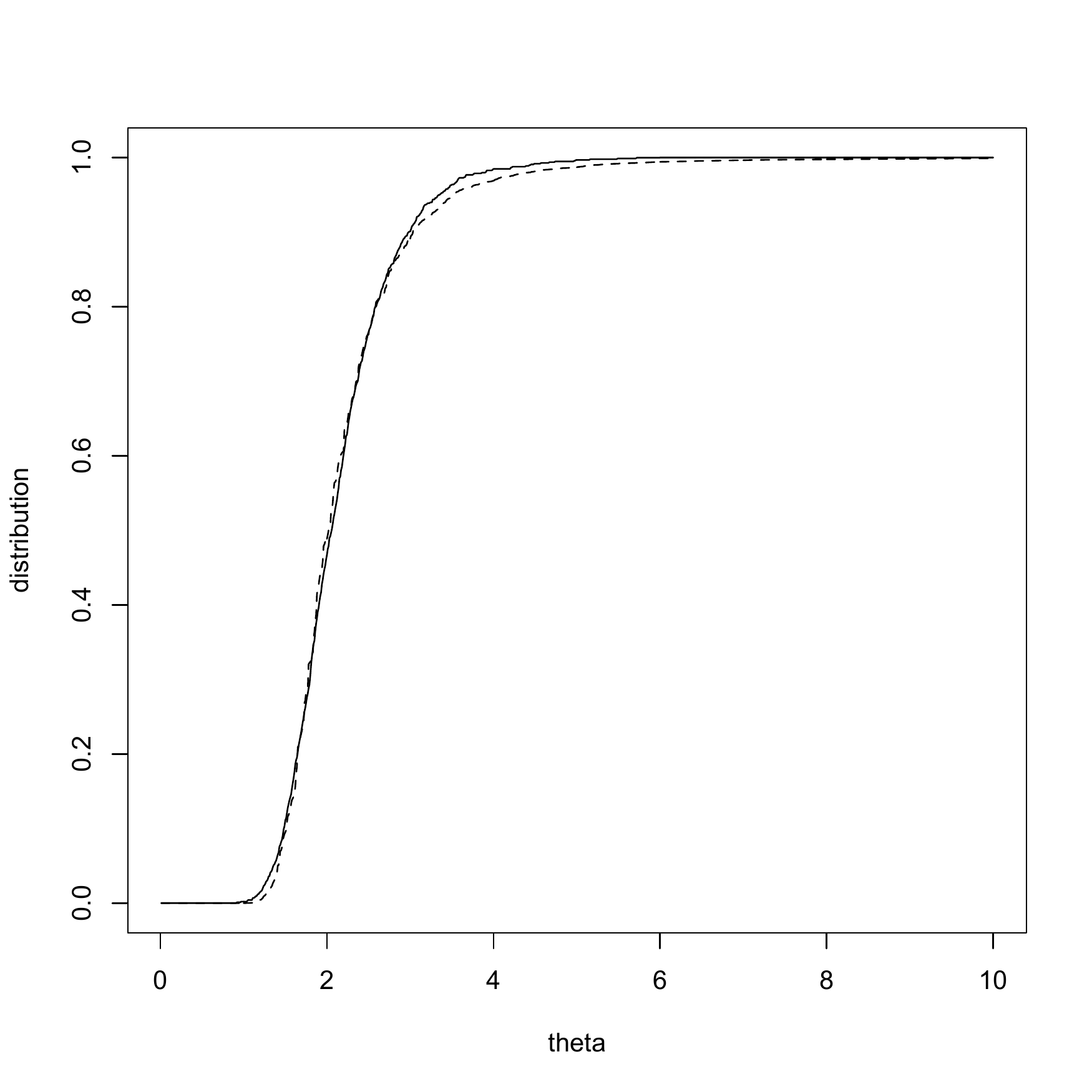}
\caption{(i) Solid line: $F^*_{\widehat{\theta}}(z)$; (ii) dashed line: $F_{\widehat{\theta}(z)}$.}
 \label{fig2}
\end{center}
\end{figure}
\end{center}

The aim here is to compare the true distribution of $\widehat{\theta}$; i.e. $F^*_{\widehat{\theta}}(z)$, based on a sample of size $n=10$ with the estimate given by (\ref{alt}). We obtain $F^*_{\widehat{\theta}}(z)$ by simulating samples of size $10$ with a true $\theta^*=2$. Repeating this multiple times and maximizing the likelihood each time yields a sample of mle's from which we construct the empirical distribution.

On the other hand, we compute $\widehat{F}_{\widehat{\theta}}(z)$ by estimating $D(z,\theta^*)$ and $V(z,\theta^*)$ arbitrarily accurately using Monte Carlo methods. The two distributions are plotted in Fig.~2; the bold line is the true distribution while the dashed line is (\ref{alt}). As can be seen, they are remarkably close for a sample of size $n=10$.

\subsection{Skew normal distribution} Here we consider the skew normal density, see \cite{azzalini}, with density function
$$f(x;\theta)=2\phi(x)\,\Phi(\theta\,x),$$
with $\theta\in\mathbb{R}$. Then $l(x,\theta)=-\log f(x;\theta)$ which is convex in $\theta$ for all $x\in\mathbb{R}$. Here we compare $T_N=\sqrt{n}\,\widehat{\theta}\,\sqrt{I(0)},$
where $I(0)=(\phi(0)/\Phi(0))^2$ is the Fisher information evaluated at $\theta=0$, with 
$$T=\sqrt{n}\,\frac{D(\widehat{\theta},0)}{\sqrt{V(\widehat{\theta},0)-D^2(\widehat{\theta},0)}}.$$
Specifically, we aim to investigate which is closer to a standard normal variable, where $\widehat{\theta}$ is obtained from a sample with $\theta=0$. 

We fix $n=15$ and generate 5000 data sets with this sample size from a standard normal distribution. This gives us 500 mle's $\widehat{\theta}$ which in turn give us 5000 values of $T_N$ and 5000 values of $T$. The first two moments of $T_N$ are $(-0.003,1.797)$ while the first two moments for $T$ are 
$(-0.064,0.799)$.  So both miss the second moment, but $T_N$ over--estimates substantially.

\begin{center}
\begin{figure}[!htbp]
\begin{center}
\includegraphics[scale=0.55]{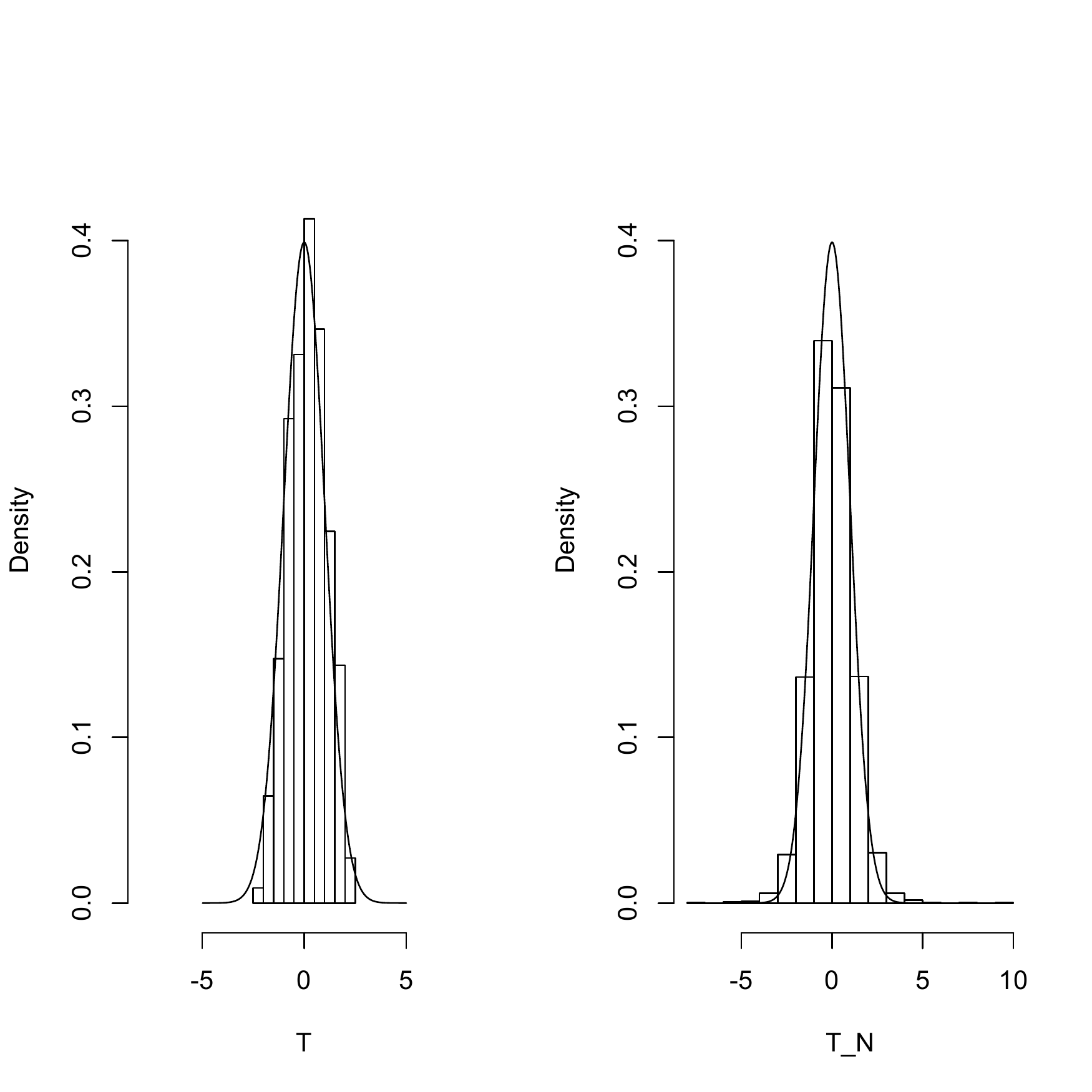}
\caption{Density estimates of $T$ and $T_N$ with $n=25$.}
 \label{fig3}
\end{center}
\end{figure}
\end{center}

When $n=25$ we plot the densities of $T$ and $T_N$, with the standard normal density overlain, in Fig.~3.
As can be seen, the $T$ still has slightly lighter tails to normal (a variance of $0.856$) but the $T_N$ remains with heavy tails (a variance of $1.379$). 
With $n=100$, we get the corresponding moment values for $T_N$ as $(-0.022,1.263)$ and for $T$ as $(-0.033,1.240)$. 

\section{Weighted likelihood bootstrap}

The weighted likelihood bootstrap was introduced in \cite{Newt} and is a way of providing approximate posterior samples. Recently the idea has had a resurgence of interest; see \cite{Newton} and \cite{Lyddon}.

The purpose of this section is to provide the exact distribution of a sample from the weighted likelihood bootstrap, in the case when 
$-\log f(x;\theta)$
is strictly convex in $\theta$ for all $x$. Here $f(x;\theta)$, with $x\in\mathbb{X}$, and $\theta\in\Theta$, is a family of density functions with $\Theta$ a one dimensional parameter space. 
Even in this case, a full anayltical study of the weighted likelihood bootstrap has not been done. Indeed, in  Section 4 of \cite{Newt} a first order approximation to a proper Bayesian procedure is all that is obtained. 

The weighted likelihood bootstrap draws a sample $\theta$ by minimizing
$$l_w(\theta)=\sum_{i=1}^n w_i\,l(x_i;\theta)$$
where $l(x;\theta)=-\log f(x;\theta)$, and the $w=(w_i)_{i=1:n}$ are from a Dirichlet distribution with all parameters set to 1;
i.e.
$$p(w)\propto {\bf 1}\left(w_i\geq 0,\,\,\sum_{i=1}^n w_i=1\right).$$
Another sample is taken by resampling $w$. In practice one can take
$$w_i=\frac{v_i}{\sum_{i=1}^n v_i}$$
where the $(v_i)$ are independent and identically distributed as standard exponential, and minimize
$$l_v(\theta)=\sum_{i=1}^n v_i\,l(x_i;\theta).$$
So note that the randomness is generated by the weights now rather than the data.

The aim is  to find $F(z)=\mbox{P}(\theta\leq z)$. This result relies partly on knowing the distribution of sums of independent exponentials; e.g.
$$S=\sum_{i=1}^n v_i\,\,\psi_i$$
where the $\psi_i>0$. See, for example, \cite{Feller}.


\subsection{Derivation of $F(z)$} The starting point is  the observation that
\begin{equation}\label{Fdist}
\mbox{P}(\theta\leq z)=\mbox{P}\left(\sum_{i=1}^n w_i\,\gamma_i(z)\geq 0\right).
\end{equation}
This follows due to the convexity of $l(x;\theta)$. Hence, we are interested in the distribution of
$$S(z)=\sum_{i=1}^n w_i\,\gamma_i(z)$$
and in particular
$$F(z)=\mbox{P}(S(z)\geq 0).$$
Since we are only interested in the probability of $S(z)$ being positive, we can represent the $(w_i)$ without their normalizing constant $\sum_{i=1:n} w_i$, and so we can take them as independent and identically distributed standard exponential random variables, $(v_i)$.

Now let us arrange $S(z)=S_1(z)-S_2(z)$ where
$$S_1(z)=\sum_{\gamma_i(z)>0}v_i\,\gamma_i(z)\quad\mbox{and}\quad S_2(z)=\sum_{\gamma_i(z)<0}v_i\,|\gamma_i(z)|$$
and
$$\gamma_i(z)=l'(x_i;z),$$
where $l'$ denotes differentiation with respect to $\theta$.
If we now present the labels so that for $i=1,\ldots,m$ it is that $\gamma_i(z)>0$, and for $i=m+1,\ldots,n$, for some $m\in\{0,\ldots,n\}$, it is that $\gamma_i(z)<0$; then define $\lambda_i(z)=1/|\gamma_i(z)|$. We can assume all the $\lambda_i$ are mutually distinct arising from the $(x_i)$ being continuous random variables.

The density function for $S_1(z)$ is
$$f_1(t;\gamma_1(z),\ldots,\gamma_m(z))=\left[\prod_{i=1}^m \lambda_i(z)\right]\,\sum_{i=1}^m q_{1i}(z)\,e^{-\lambda_i(z) t},$$
where, for $i=1,\ldots,m$,
$$q_{1i}(z)=\prod_{k=1:m,\,k\ne i}\frac{1}{\lambda_k(z)-\lambda_i(z)}.$$
See \cite{Feller}. Likewise, the density function for $S_2(z)$ is
$$f_2(t;\gamma_{m+1}(z),\ldots,\gamma_n(z))=\left[\prod_{i=m+1}^n \lambda_i(z)\right]\,\sum_{i=m+1}^n q_{2i}(z)\,e^{-\lambda_i(z) t},$$
where, for $i=m+1,\ldots,n$,
$$q_{2i}(z)=\prod_{k=m+1:n,\,k\ne i}\frac{1}{\lambda_k(z)-\lambda_i(z)}.$$
Hence, it is now easy to see that
$$F(z)=\int_0^\infty \bar{F}_1\big(t;\gamma_1(z),\ldots,\gamma_m(z)\big)\,f_2(t;\gamma_{m+1}(z),\ldots,\gamma_n(z))\,d t,$$
where $\bar{F}$ represents the survival function corresponding to density $f$.
The integration is straightforward leading to
\begin{equation}\label{Fexact}
F(z)=\left(\prod_{i=1}^n \lambda_i(z)\right)\,\sum_{l=1}^m\sum_{j=m+1}^n \frac{q_{1l}(z)\,q_{2j}(z)}{\lambda_l(z)\,(\lambda_l(z)+\lambda_j(z))}.
\end{equation}
While a complicated function of $z$, it is quite easy to compute numerically.

\begin{center}
\begin{figure}[!htbp]
\begin{center}
\includegraphics[scale=0.48]{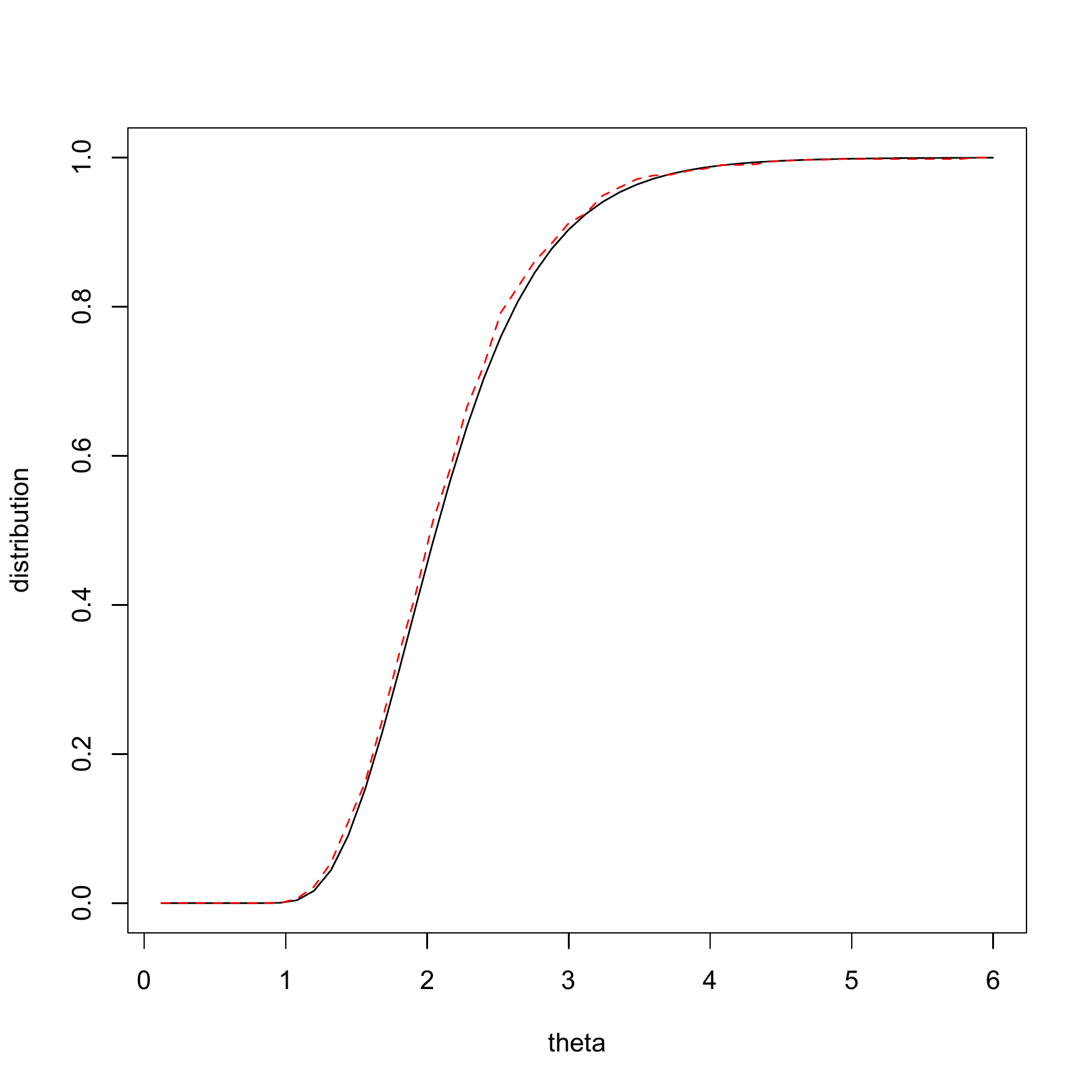}
\caption{Exact (bold line) and estimated (dashed red line) weighted likelihood bootstrap posterior distribution for the beta model}
\label{fig3}
\end{center}
\end{figure}
\end{center}

In Fig.~\ref{fig3} we present the exact weighted likelihood bootstrap distribution for the model $f(x;\theta)=\theta\,x^{\theta-1}$, with $\theta>0$ and $0<x<1$. In this case $l'(x;\theta)=-1/\theta-\log x$, and we took $n=10$ samples from a beta$(2,1)$ distribution. The bold line is the exact distribution of the weighted likelihood posterior for $\theta$ and the red dashed line is the approximate distribution obtained from 1000 samples from the weighted likelihood bootstrap.

\begin{center}
\begin{figure}[!htbp]
\begin{center}
\includegraphics[scale=0.5]{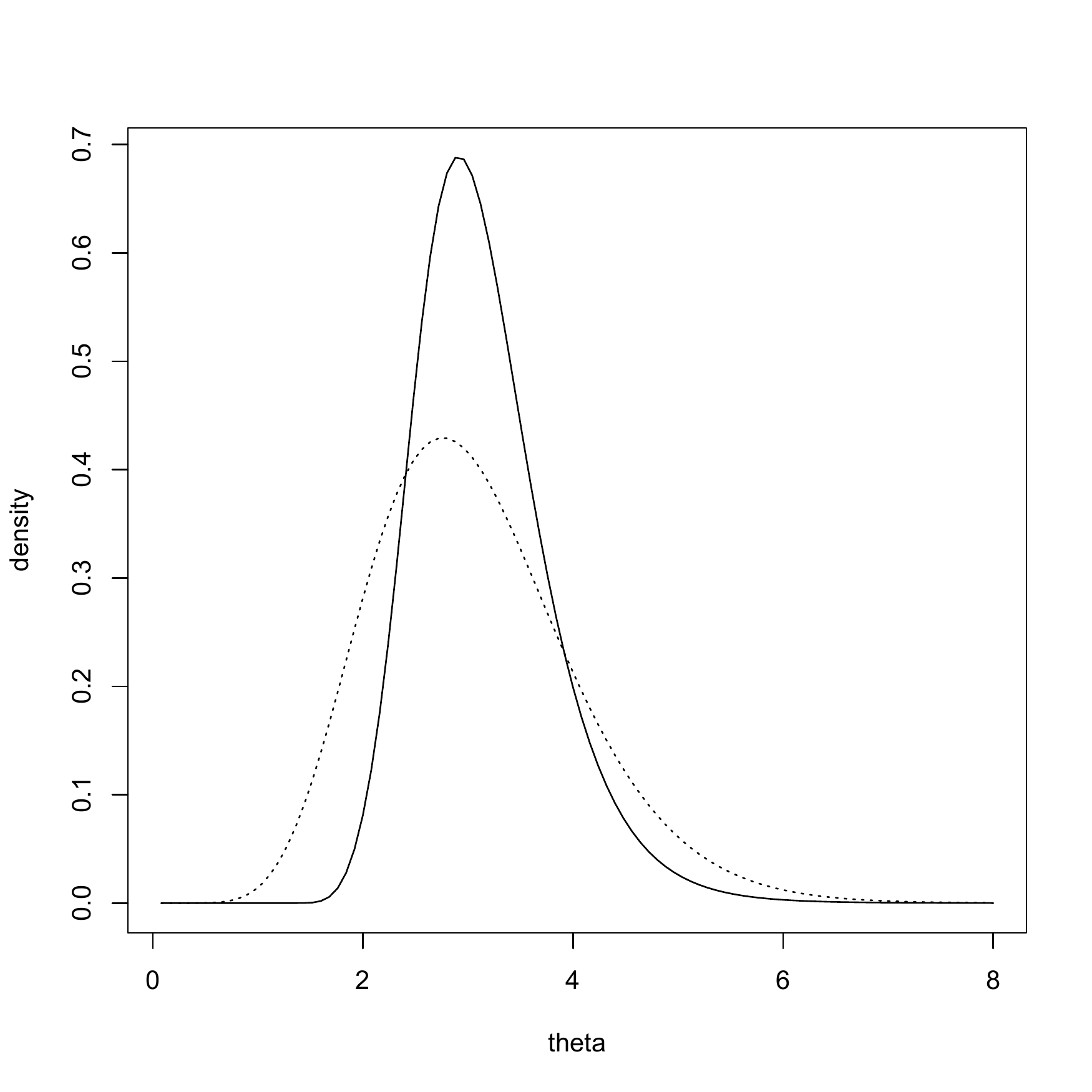}
\caption{Exact (bold line) weighted likelihood bootstrap posterior density and exact posterior with Jeffrey's prior (dashed line) for the exponential  model}
\label{fig4}
\end{center}
\end{figure}
\end{center}

In Fig.~\ref{fig4} we present the exact weighted likelihood bootstrap posterior density, obtained by the numerical differentiation of (\ref{Fexact}, compared with the posterior density  using Jeffrey's prior. The model used is exponential; i.e. $f(x;\theta)=\theta\,\exp(-x\theta)$ and the $n=10$ data points were taken from this model with a true $\theta=1/3$. 

In this example it is clear that the weighted likelihood bootstrap has less posterior variance compared to that provided by the Jeffrey's prior.

\subsection{Asymptotic approximation} On inspection of (\ref{Fdist}) one can see that a normal type approximation (though not a normal distribution for $F$) is going to be provided by
\begin{equation}\label{Fapprox}
\widehat{F}(z)= \Phi\left(\frac{\sum_{i=1}^n \gamma_i(z)}{\sqrt{\sum_{i=1}^n \gamma_i^2(z)}}\right),
\end{equation}
where $\Phi$ denotes the standard normal cumulative function.  This follows from standard asymptotic theory; namely that, for large $n$,
$$S_n=n^{-1}\sum_{i=1}^n v_i\,\gamma_i(z)$$
will be approximately normal with mean and variance given by
$$\mbox{E}\,S_n=n^{-1}\sum_{i=1}^n\gamma_i(z)\quad\mbox{and}\quad \mbox{Var}\,S_n=n^{-2}\sum_{i=1}^n \gamma_i^2(z),$$
respectively. Then (\ref{Fapprox}) follows since $\mbox{Pr}(\mbox{Z}(\mu,\sigma^2)\geq 0)=\Phi(\mu/\sigma)$, where $\mbox{Z}(\mu,\sigma^2)$ denotes a normal random variable with mean $\mu$ and variance $\sigma^2$. 

We can develop the asymptotic approximation further, relying on
$$0=n^{-1}\sum_{i=1}^n l'(x_i;\widehat{\theta})=n^{-1}\sum_{i=1}^n l'(x_i;\theta)+(\widehat{\theta}-\theta)\,n^{-1}\sum_{i=1}^n l''(x_i;\theta)+o(|\widehat{\theta}-\theta|),$$
for small $|\widehat{\theta}-\theta|$, where $\widehat{\theta}$ is the maximum likelihood estimator. 
Given that
$$n^{-1}\sum_{i=1}^n l''(x_i;\theta)\quad\mbox{and}\quad n^{-1}\sum_{i=1}^n \big(l'(x_i;\theta)\big)^2$$
are asymptotically equivalent, both approximating the Fisher information, $I(\theta)$, 
we obtain the asymptotic equivalence between
$$\frac{\sum_{i=1}^n \gamma_i(z)}{\sqrt{\sum_{i=1}^n \gamma_i^2(z)}}\quad\mbox{and}\quad  \sqrt{n}\,(z-\widehat{\theta})\,\sqrt{I(z)}.$$
Hence, a further asymptotic approximation to (\ref{Fapprox}) is given by
\begin{equation}\label{Ftilde}
\widetilde{F}(z)=\Phi(T_n(z))
\end{equation}
where
$$T_n(z)=\sqrt{n}\,\left(z-\widehat{\theta}\right)\,\sqrt{I(z)}.$$
It is possibe to see  (\ref{Fapprox}) as a Bayesian probability matching type procedure. See, for example, \cite{Datta}. In particular, in Section 3 of \cite{Ghosh}, the authors consider 
$$T_n(\theta)=\sqrt{n}\,(\theta-\widehat{\theta})\,\sqrt{I(\theta)}.$$
The probability matching idea is to treat $T_n$ in two ways; first as a random object induced by the random sample with $\theta$ as the fixed true value, and, second, as random, with the data now fixed and the randomness induced by a posterior distribution on $\theta$, having found a suitable prior $\pi(\theta)$, to ensure
$$\mbox{E}\left[ \mbox{P}^{\pi}(T_n\leq z|x_1,\ldots,x_n)|\theta\right]=\mbox{P}(T_n\leq z|\theta)+o(1/n).$$ 
The former interpretation has $T_n$ as asymptotically standard normal, a well known result.

The need to provide such a matching via a prior to posterior procedure has recently been challenged; see for example \cite{Belitser}. Accepting the idea that one can directly construct a posterior, one can obtain  (\ref{Ftilde}) directly as a posterior which is not based on any prior. In this scenario, an asymptotic motivated probability matching ``posterior'', asymptotically equivalent to the weighted likelihood bootstrap, is provided by samples $\theta$, whereby $\theta$ solves
\begin{equation}\label{probmatch}
\sqrt{n}(\theta-\widehat{\theta})\sqrt{I(\theta)}=z,
\end{equation}
where $z$ is a  standard normal random variable. 

Here we compare samples from (\ref{probmatch}) with those from a weighted likelihood bootstrap with $n=100$ data points and model $f(x;\theta)=\exp(\theta-xe^\theta)$, with $\theta\in (-\infty,+\infty)$. We took the data with true parameter as $\log 3$ and took 100 samples from (\ref{probmatch}) and the weighted likelihood bootstrap. The two empirical distributions of the samples are presented in Fig.~\ref{fig2}. The solid bold line is for the weighted likelihood bootstrap and the dashed line for the (\ref{probmatch}) samples.

\begin{center}
\begin{figure}[!htbp]
\begin{center}
\includegraphics[scale=0.5]{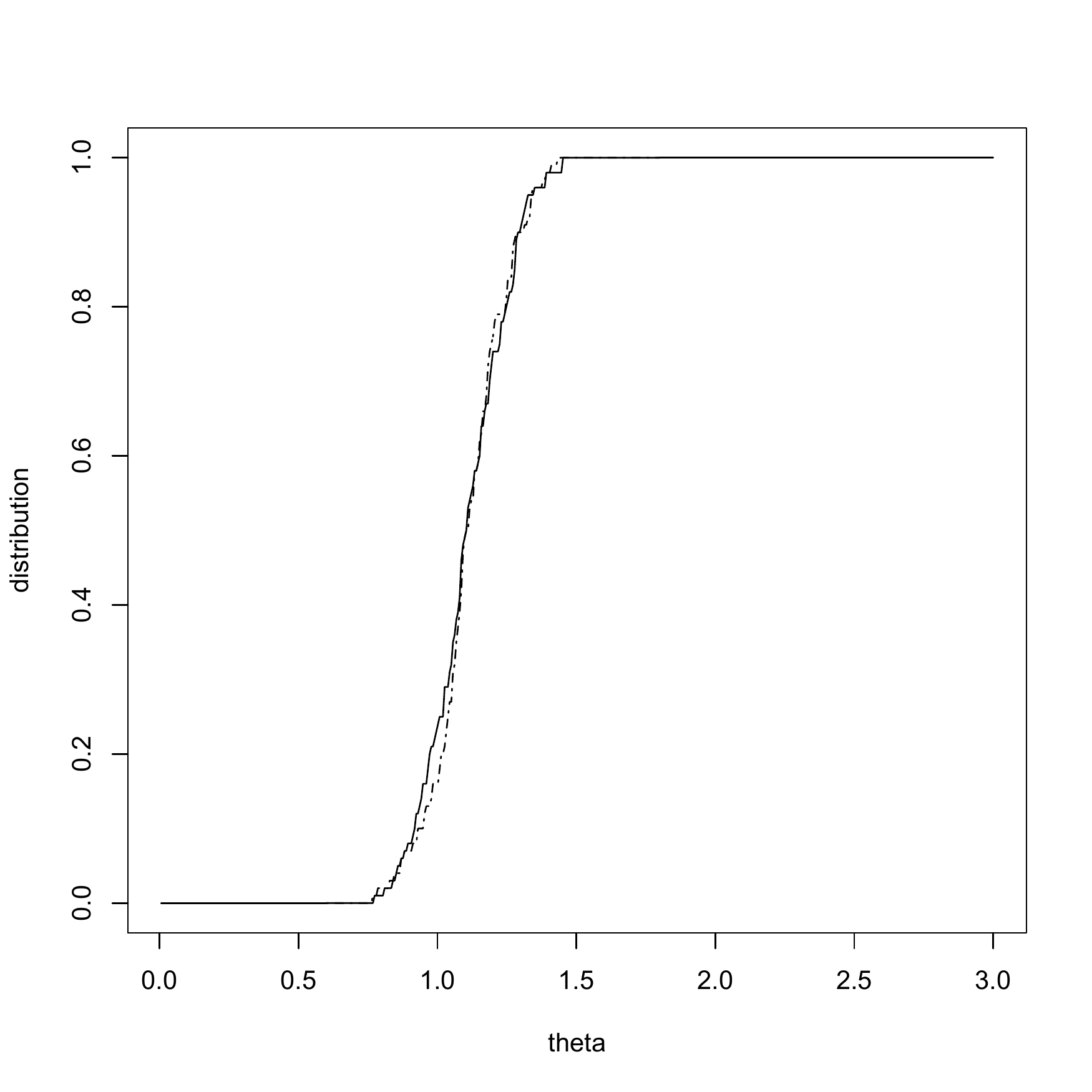}
\caption{Empirical distributions corresponding to samples from weighted likelihood bootstrap (solid) and probability matching posterior (dashed) for exponential model}
\label{fig2}
\end{center}
\end{figure}
\end{center}

\section{Discussion} If $l(x;\theta)$ is strictly convex in $\theta$ for all $x$ then we can obtain an accurate estimate of the distribution of the maximum likelihood estimate using only $l'(x;\theta)$. For the multivariate case; i.e. $\Theta\subset \mathbb{R}^d$, it is not easy in general to find
\begin{equation}\label{mult}
F_{\widehat{\theta}}(z_1,\ldots,z_d)=\P\left(\sum_{i=1}^n\frac{\partial}{\partial \theta_1}l(x_i;z_1)\geq 0, \ldots, \sum_{i=1}^n\frac{\partial}{\partial \theta_d}l(x_i;z_d) \geq 0 \right).
\end{equation}
Using a multivariate normal approximation, say $\mbox{MVN}_d(\mu(z),\Sigma(z))$, to the vector $\psi(z)$,  where
$$\psi_j(z)=\frac{1}{n}\sum_{i=1}^n\frac{\partial}{\partial \theta_j}l(x_i;z_j),$$
we would have
$F_{\widehat{\theta}}(z)\approx \P(Y(z)\geq 0) $ with $Y(z)\sim \mbox{MVN}_d(\mu(z),\Sigma(z))$. Here
$$\mu_j(z)=\int \frac{\partial}{\partial\theta_j}l(x;z_j)\,f(x;\theta^*)\,\d x,\quad\Sigma_{jk}=
\int\frac{\partial}{\partial\theta_j} l(x;z_j)\,\frac{\partial}{\partial\theta_k} l(x;z_k)\,f(x;\theta^*)\,d x.$$
Approximating $\P(Y(z)\geq 0)$ in multidimensions has been considered, for example, by \cite{cox}, who could only find adequate approximations up to 3 dimensions.

An approximate sampling strategy from (\ref{mult}) would involve the parametric bootstrap; see for example \cite{efron}. 
To sample $z$ from (\ref{mult}) approximately, take a sample $\widetilde{x}=(\widetilde{x}_1,\ldots,\widetilde{x}_n)$ from $f(\cdot;\widehat{\theta})$ and take $z$ as the mle with data $\widetilde{x}$; i.e.
take 
$$z=\arg\min_\theta\sum_{1\leq i\leq n} l(\widetilde{x}_i;\theta)$$
as approximately coming from (\ref{mult}).

On the other hand, for the multivariate weighted likelihood bootstrap we can use a sequence of conditional densities. So now assume that $l(x;\theta)$, with $\theta\in\Theta\subset\mathbb{R}^d$, is such that
$\partial^2 l(x;\theta)/\partial\theta_j^2\geq 0$ for all $j=1,\ldots,d$ and all $x$. Then
$$\mbox{P}(\theta_j\leq z_j|\theta_{-j}=z_{-j})=\mbox{P}\left(\sum_{i=1}^n v_i\,\partial l(x_i;z)/\partial\theta_j\geq 0\right),$$
where $z=(z_1,\ldots,z_d)=(z_j,z_{-j})$. Hence, we can find easily the conditional density equivalent of (\ref{Fexact}); i.e. $F(z_j|z_{-j})$ for each $j\in\{1,\ldots,d\}$.


\end{document}